\documentclass[twoside,10pt]{article}
\usepackage{amssymb}
\usepackage{amsmath,amssymb,amsthm}

\begin{document}

\setlength{\textwidth}{126mm} \setlength{\textheight}{180mm}
\setlength{\parindent}{0mm} \setlength{\parskip}{2pt plus 2pt}

\frenchspacing

\pagestyle{myheadings}

\markboth{Dimitar Mekerov}{On Riemannian almost product manifolds with nonintegrable structure}

\newtheorem{thm}{Theorem}[section]
\newtheorem{lem}[thm]{Lemma}
\newtheorem{prop}[thm]{Proposition}
\newtheorem{cor}[thm]{Corollary}
\newtheorem{probl}[thm]{Problem}

\newtheorem{defn}{Definition}[section]
\newtheorem{rem}{Remark}[section]
\newtheorem{exa}{Example}



\newcommand{\X}{\mathfrak{X}}
\newcommand{\B}{\mathcal{B}}
\newcommand{\s}{\mathfrak{S}}
\newcommand{\g}{\mathfrak{g}}
\newcommand{\W}{\mathcal{W}}
\newcommand{\Lgr}{\mathrm{L}}
\newcommand{\dd}{\mathrm{d}}

\newcommand{\pd}{\partial}
\newcommand{\ddx}{\frac{\pd}{\pd x^i}}
\newcommand{\ddy}{\frac{\pd}{\pd y^i}}
\newcommand{\ddu}{\frac{\pd}{\pd u^i}}
\newcommand{\ddv}{\frac{\pd}{\pd v^i}}

\newcommand{\diag}{\mathrm{diag}}
\newcommand{\End}{\mathrm{End}}
\newcommand{\im}{\mathrm{Im}}
\newcommand{\id}{\mathrm{id}}

\newcommand{\ie}{i.e.}
\newfont{\w}{msbm9 scaled\magstep1}
\def\R{\mbox{\w R}}
\newcommand{\norm}[1]{\left\Vert#1\right\Vert ^2}
\newcommand{\nN}{\norm{N}}
\newcommand{\nP}{\norm{\nabla P}}
\newcommand{\tr}{{\rm tr}}

\newcommand{\nJ}[1]{\norm{\nabla J_{#1}}}
\newcommand{\thmref}[1]{Theorem~\ref{#1}}
\newcommand{\propref}[1]{Proposition~\ref{#1}}
\newcommand{\secref}[1]{\S\ref{#1}}
\newcommand{\lemref}[1]{Lemma~\ref{#1}}
\newcommand{\dfnref}[1]{Definition~\ref{#1}}

\frenchspacing


\title{On Riemannian almost product manifolds \\ with nonintegrable structure}

\author{Dimitar Mekerov}

\maketitle

{\small
{\it Abstract.} The class of the Riemannian almost product
manifolds with nonintegrable structure is considered. Some
identities for curvature tensor as certain invariant tensors and
quantities are obtained.

{\it Mathematics Subject Classification (2000):} 53C15, 53C50   \\
{\it Key words:} almost product manifold, Riemannian metric,
nonintegrable structure}


\section*{Introduction}

The systematic development of the theory of Riemannian almost
product manifolds was started by K. Yano \cite{Yano}. In
\cite{Nav} A.~M.~Naveira gives a classification of these manifolds
with respect to the covariant differentiation of the almost product
structure.

Such manifolds with zero trace of the almost product structure are
considered in \cite{Mih}, \cite{St:prop}--\cite{StGrMe:const}.
Moreover, a classification is presented in \cite{StGr:connect},
having in mind the results in \cite{Nav}. In the classification in
\cite{StGr:connect} the basic class $\W_3$ is only the class with
nonintegrable structure.

In the present work the problems in the differential geometry of
the manifolds of the class $\W_3$ are mainly considered.


\section{Preliminaries}

Let $(M,P,g)$ be a Riemannian almost product manifold, \ie{} a
differentiable manifold $M$ with a tensor field $P$ of type
$(1,1)$ and a Riemannian metric $g$ such that
\begin{equation}\label{Pg}
    P^2x=x,\quad g(Px,Py)=g(x,y)
\end{equation}
for arbitrary $x$, $y$ of the algebra $\X(M)$ of the smooth vector
fields on $M$. Obviously $g(Px,y)=g(x,Py)$.

Further $x,y,z,w$ will stand for arbitrary elements of $\X(M)$.

In this work we consider Riemannian almost product manifolds with
$\tr{P}=0$. In this case $(M,P,g)$ is an even-dimensional
manifold.

If $\dim{M}=2n$ then the associated metric $\tilde{g}$ of $g$, determined by $\tilde{g}(x,y)=g(x,Py)$,
is an indefinite metric of signature $(n,n)$.
Since $\tilde{g}(Px,Py)=\tilde{g}(x,y)$, the manifold $(M,P,\tilde{g})$ is a pseudo-Riemannian almost
product manifold.
We say that $(M,P,\tilde{g})$ is an \emph{associated manifold} of $(M,P,g)$.

The classification from \cite{StGr:connect} of Riemannian almost
product manifolds is made with respect to the tensor field $F$ of
type (0,3), defined by
\begin{equation}\label{2}
F(x,y,z)=g\left(\left(\nabla_x P\right)y,z\right),
\end{equation}
where $\nabla$ is the Levi-Civita connection of $g$. The tensor $F$ has the following properties:
\begin{equation}\label{3}
    F(x,y,z)=F(x,z,y)=-F(x,Py,Pz),\quad F(x,y,Pz)=-F(x,Py,z).
\end{equation}
The only class of Riemannian almost product manifolds with
nonintegrable structure is the basic class $\W_3$ determined by
the condition
\begin{equation}\label{sigma}
    \mathop{\s}_{x,y,z} F(x,y,z)=0,
\end{equation}
where $\mathop{\s}_{x,y,z}$ is the cyclic sum by $x, y, z$.

Further manifolds of the class $\W_3$ we call \emph{Riemannian
$\W_3$-manifolds}.

The condition \eqref{sigma} is equivalent to
\begin{equation}\label{sigmaP}
        \mathop{\s}_{x,y,z} F(Px,y,z)=0.
\end{equation}

In \cite{StGr:connect} the symmetric tensor field $\bar{N}$ is defined by
\begin{equation}\label{6}
    \bar{N}(x,y)=\left(\nabla_x P\right)Py+\left(\nabla_{Px} P\right)y+\left(\nabla_y P\right)Px+\left(\nabla_{Py} P\right)x
\end{equation}
and it has the properties
\[
    \bar{N}(Px,Py)=\bar{N}(x,y), \quad \bar{N}(Px,y)=\bar{N}(x,Py), \quad \bar{N}(x,Py)=-P\bar{N}(x,y).
\]
It is proved that the condition \eqref{sigma} is equivalent to
$\bar{N}(x,y)=0$.

The class $\W_0$, defined by the condition $F(x,y,z)=0$, is
contained in the other classes. This is the class of so-called
\emph{Riemannian $P$-manifolds}, i.e. differentiable
even-dimensional manifolds $(M,P,g)$ with Riemannian metric $g$
and structure $P$, such that $g(Px,Py)=g(x,y)$, $P^2=\id$,
$\tr{P}=0$, $\nabla P=0$. Therefore the class $\W_0$ is an
analogue of the class of K\"ahlerian manifolds in the geometry of
almost Hermitian manifolds.

The following property of the covariant derivation of $F$ is valid
\begin{equation}\label{8}
    \left(\nabla_x F\right)(y,z,Pw)+\left(\nabla_x F\right)(y,Pz,w)=A(x,y,z,w),
\end{equation}
where
\begin{equation}
    A(x,y,z,w)=-g\bigl(\left(\nabla_x P\right)z,\left(\nabla_y P\right)w\bigr)-g\bigl(\left(\nabla_y P\right)z,\left(\nabla_x P\right)w\bigr).
\end{equation}

As it is known the curvature tensor field $R$ of a Riemannian
manifold with metric $g$ is determined by
\[
    R(x,y)z=\nabla_x \nabla_y z - \nabla_y \nabla_x z -
    \nabla_{[x,y]}z
\]
and the corresponding tensor field of type $(0,4)$ is defined as
follows
\[
    R(x,y,z,w)=g(R(x,y)z,w).
\]

Let $(M,P,g)$ be a Riemannian almost product manifold and
$\{e_i\}$ be a basis of the tangent space $T_pM$ at a point $p\in
M$. Let the components of the inverse matrix of $g$ with respect
to $\{e_i\}$ be $g^{ij}$. If $\rho$ and $\tau$ are the Ricci
tensor and the scalar curvature, then $\rho^*$ and $\tau^*$,
defined by $\rho^*(y,z)=g^{ij}R(e_i,y,z,Pe_j)$ and
$\tau^*=g^{ij}\rho^*(e_i,e_j)$, are called an \emph{associated
Ricci tensor} and an \emph{associated scalar curvature},
respectively. We denote
$\tau^{**}=g^{ij}g^{kl}R(e_i,e_k,Pe_l,Pe_j)$.

The Lie form $\theta$ associated to $F$ is defined by $\theta
(z)=g^{ij}F(e_i,e_j,z)$. For a Riemannian $\W_3$-manifold it is
valid $\theta(z)=0$.

The square norm of $\nabla P$ is defined by
$\nP=g^{ij}g^{kl}g\left(\left(\nabla_{e_i}P\right)e_k,\left(\nabla_{e_j}P\right)e_l\right)$.


\section{A basic identity for the curvature tensor of Riemannian $\W_3$-manifolds}

Let $(M,P,g)$ be a Riemannian $\W_3$-manifold. According to \eqref{sigmaP} we have
\begin{equation}\label{sigmaP2}
            \mathop{\s}_{y,z,w} F(Py,z,w)=0.
\end{equation}
Then the covariant differentiations of \eqref{sigmaP2} and
\eqref{sigma} imply, respectively, the following equations
\begin{equation}\label{11}
    \mathop{\s}_{y,z,w}\left(\nabla_x F\right)(Py,z,w)=
    \mathop{\s}_{y,z,w}g\bigl(\left(\nabla_x P\right)y,\left(\nabla_z P\right)w+\left(\nabla_w P\right)z\bigr),
\end{equation}
\begin{equation}\label{12}
    \mathop{\s}_{y,z,w}\left(\nabla_x F\right)(y,z,w)=0.
\end{equation}

Having in mind the Ricci identity
\[
\left(\nabla_x \nabla_y P\right)z-\left(\nabla_y \nabla_x P\right)z=R(x,y)Pz-PR(x,y)z
\]
and the properties \eqref{3} it follows immediately
\begin{equation}\label{13}
    \left(\nabla_x F\right)(y,z,w)-\left(\nabla_y F\right)(x,z,w)=R(x,y,Pz,w)-R(x,y,z,Pw).
\end{equation}
Using \eqref{13} we get that the left-hand side of \eqref{11} has the form
\begin{equation}\label{14}
    \mathop{\s}_{y,z,w}\left(\nabla_{Py} F\right)(x,z,w)+\mathop{\s}_{y,z,w}\left\{R(x,Py,Pz,w)-R(x,Py,z,Pw)\right\}.
\end{equation}

We apply \eqref{12} to each of the addends of the first cyclic sum
in \eqref{14} and then we apply \eqref{11}, \eqref{12}, \eqref{13}
in the obtained expression. After that we use the properties of
$R$ and we finally get the following identity
\begin{equation}\label{15}
\begin{split}
    \mathop{\s}_{x,y,z}&\left\{R(x,Py,Pz,w)-R(x,Py,z,Pw)\right. \\
    \phantom{\mathop{\s}_{x,y,z}}
    &\left. +R(Px,y,z,Pw)-R(Px,y,Pz,w)\right\}\\
    =&\mathop{\s}_{x,y,z} g\bigl(\left(\nabla_x P\right)y+\left(\nabla_y P\right)x,\left(\nabla_z P\right)w+\left(\nabla_w P\right)z\bigr).
\end{split}
\end{equation}

In this way we prove the following
\begin{thm}\label{thm:15}
If $(M,P,g)$ is a Riemannian $\W_3$-manifold, then the curvature tensor $R$ has the property \eqref{15}.
\end{thm}

Having in mind the last theorem we prove the following
\begin{cor}\label{thm:15-trace}
If $(M,P,g)$ is a Riemannian $\W_3$-manifold, then the following
equations are valid
\begin{equation}\label{rho}
\begin{split}
    &\rho(y,z)+\rho(Py,Pz)-\rho^*(Py,z)-\rho^*(y,Pz)=\\[4pt]
    &=g^{ij}g\bigl(\left(\nabla_{e_i} P\right)y+\left(\nabla_y P\right)e_i,\left(\nabla_z P\right)e_j
    +\left(\nabla_{e_j} P\right)z\bigr),
\end{split}
\end{equation}
\begin{equation}\label{nP}
    \nP=-2g^{ij}g^{kl}g\left(\left(\nabla_{e_i}P\right)e_k,\left(\nabla_{e_l}P\right)e_j\right),
\end{equation}
\begin{equation}\label{nP2}
    \nP=2\left(\tau-\tau^{**}\right).
\end{equation}
\end{cor}
\begin{proof}
It is satisfied \eqref{15} for a Riemannian $\W_3$-manifold hence we obtain \eqref{rho}. Having in
mind \eqref{sigma} and the definition of $\nP$, we have
\eqref{nP}. The equations \eqref{rho} and \eqref{nP} imply
\eqref{nP2}.
\end{proof}

\begin{rem}
If $(M,P,g)$ is a Riemannian $\W_3$-manifold with $\dim M\geq 4$
and a K\"ahler curvature tensor $R$, i.e.
$R(x,y,Pz,Pw)=R(x,y,z,w)$, then $\tau^{**}=\tau$ and therefore
$\nP=0$.
\end{rem}


\section{Invariant bisectional curvature}

Let $(M,P,g)$ be a Riemannian almost product manifold and let
$\alpha$ be a 2-plane in the tangent space at an arbitrary point
of $M$. Then, it is known, the \emph{sectional curvature} of
$\alpha$ is defined by the following equation
\[
    k(x,y)=\frac{R(x,y,y,x)}{g(x,x)g(y,y)-g^2(x,y)},
\]
where $(x,y)$ is an arbitrary basis of $\alpha$.
If $x$ is a noneigenvector of $P$, then the 2-plane $\alpha=(x,Px)$ is called an \emph{invariant 2-plane},
because $P\alpha=\alpha$.

Let us define the quantity
\begin{equation}\label{16}
    h(x,y)=\frac{R(x,Px,y,Py)}{\sqrt{g^2(x,x)-g^2(x,Px)}\sqrt{g^2(y,y)-g^2(y,Py)}},
\end{equation}
where $x$, $y$ are noneigenvectors of $P$.

It is known from \cite{StGrMe:invar} that an orthonormal adapted basis exists in every invariant 2-plane.
Let this basis be $(e_1,Pe_1)$ for $\alpha_1=(x,Px)$.
Then we have $x=\lambda e_1+\mu Pe_1$, $Px=\lambda Pe_1+\mu e_1$ ($\lambda, \mu \in \R$),
hence $g(x,x)=\lambda^2+\mu^2$, $g(x,Px)=2\lambda\mu$. Therefore  $g^2(x,x)-g^2(x,Px)=
\left(\lambda^2-\mu^2\right)^2\geq 0$. The assumption $\lambda^2=\mu^2$ (i.e. $\lambda=\pm\mu$)
implies $x=\pm Px$, which is a contradiction to the condition $x$ to be a noneigenvector of $P$.
Hence we get $g^2(x,x)-g^2(x,Px)>0$ and therefore the quantity $h(x,y)$ is defined correctly in \eqref{16}.

\begin{thm}
The quantity $h(x,y)$ defined by \eqref{16} for a Riemannian
almost product manifold $(M,P,g)$ depends on the invariant
2-planes $\alpha_1=(x,Px)$ and  $\alpha_2=(y,Py)$ only.
\end{thm}
\begin{proof}
 Let $z\in \alpha_1$, $w\in \alpha_2$ be noneigenvectors of $P$. Then $z=\lambda_1 x+\mu_1 Px$,
 $w=\lambda_2 y+\mu_2 Py$ ($\lambda_1, \lambda_2, \mu_1, \mu_2 \in \R$), from where we have
\[
\begin{array}{c}
g^2(z,z)-g^2(z,Pz)=\left(\lambda_1^2-\mu_1^2\right)\left(g^2(x,x)-g^2(x,Px)\right),\\[4pt]
g^2(w,w)-g^2(w,Pw)=\left(\lambda_2^2-\mu_2^2\right)\left(g^2(y,y)-g^2(y,Py)\right),\\[4pt]
R(z,Pz,w,Pw)=\left(\lambda_1^2-\mu_1^2\right)\left(\lambda_2^2-\mu_2^2\right)R(x,Px,y,Py).\\[4pt]
\end{array}
\]
Therefore   $h(x,y)=h(z,w)$.
\end{proof}

The quantity $h(x,y)$ defined by \eqref{16} we call an
\emph{invariant bisectional curvature} of invariant 2-planes
$\alpha_1=(x,Px)$ and $\alpha_2=(y,Py)$. In particular, if $x=y$,
the 2-planes $\alpha_1$ and $\alpha_2$ coincide and then
$h(x,x)=k(x,Px)$, i.e. $h(x,x)$ is the Riemannian sectional
curvature of the invariant 2-plane $(x,Px)$.

\begin{thm}\label{thm:17}
A Riemannian almost product manifold $(M,P,g)$ has a zero invariant bisectional curvature if
and only if the following condition is satisfied
\begin{equation}\label{17}
    R(x,Py,Pz,w)-R(x,Py,z,Pw)+R(Px,y,z,Pw)-R(Px,y,Pz,w)=0.
\end{equation}
\end{thm}
\begin{proof}
Let $(M,P,g)$ be a Riemannian almost product manifold with $h(x,y)=0$. Then $R(x,Px,y,Py)=0$.
In the last equation let $x+z$ and $y+w$ stand for $x$ and $y$, respectively. Hence we obtain
\[
    R(x,Pz,y,Pw)+R(x,Pz,w,Py)+R(z,Px,y,Pw)+R(z,Px,w,Py)=0.
\]
Next we substitute $y\leftrightarrow z$ in the last equation and
we get \eqref{17}.

Now, let \eqref{17} be valid. There let $x$ and $z$ stand for $y$ and $w$, respectively.
Then we obtain $R(x,Px,z,Pz)=0$, i.e. the manifold has a zero invariant bisectional curvature.
\end{proof}

\thmref{thm:15}, \thmref{thm:15-trace} and \thmref{thm:17} imply
immediately the following
\begin{cor}
If a Riemannian $\W_3$-manifold $(M,P,g)$ has a zero invariant
bisectional curvature then the following conditions are valid:
\begin{equation}\label{i}
        \mathop{\s}_{x,y,z} g\bigl(\left(\nabla_x P\right)y+\left(\nabla_y P\right)x,
            \left(\nabla_z P\right)w+\left(\nabla_w P\right)z\bigr);
\end{equation}
\begin{equation}\label{ii}
    \tau^{**}=\tau;
\end{equation}
\begin{equation}\label{iii}
        \nP=0.
\end{equation}
\end{cor}


\section{An associated pseudo-Riemannian almost \\ product manifold}

Let $(M,P,g)$ be a Riemannian almost product manifold and
$(M,P,\tilde{g})$ be the associated pseudo-Riemannian almost
product manifold. We denote the Levi-Civita connection of
$\tilde{g}$ by $\tilde{\nabla}$. In \cite{StGr:connect} the tensor
field $\Phi$ of type (1,2) is defined by
\[
\Phi(x,y)=\tilde{\nabla}_x y- \nabla_x y
\]
and the corresponding tensor field of type (0,3) by
\begin{equation}\label{22}
\Phi(x,y,z)=g\left(\tilde{\nabla}_x y- \nabla_x y,z\right).
\end{equation}

It is proved that
\begin{equation}\label{23}
\Phi(x,y,z)=\frac{1}{2}\bigl(F(x,y,Pz)+F(y,Pz,x)-F(z,x,y)\bigr).
\end{equation}

We have a classification of pseudo-Riemannian almost product
manifolds \\ $(M,P,\tilde{g})$ with respect to the tensor field
$\tilde{F}$ of type (0,3), defined by $\tilde{F}(x,y,z)=$\\
$\tilde{g}\left(\left(\tilde\nabla_x P\right) y,z\right)$. This
classification is analogous to the one of Riemannian almost
product manifolds $(M,P,g)$ with respect to $F$. The class defined
by \\ $\mathop{\s}_{x,y,z} \tilde F(x,y,z)=0$ is the only class of
nonintegrable structure $P$. This class will be denoted by $\W_3$,
too, and its manifolds will be called \emph{pseudo-Riemannian
$\W_3$-manifolds}. The class $\W_0$, defined by
$\tilde{F}(x,y,z)=0$ is contained in the other classes. This is
the class of \emph{pseudo-Riemannian $P$-manifolds}, i.e.
differentiable $2n$-dimensional manifolds $(M,P,\tilde{g})$ with
pseudo-Riemannian metric $\tilde{g}$ of signature $(n,n)$ and a
structure $P$ such that
\[
\tilde{g}(Px,Py)=\tilde{g}(x,y),\quad P^2=\id,\quad \tr{P}=0,\quad
\tilde\nabla{P}=0.
\]

\begin{thm}
A Riemannian almost product manifold $(M,P,g)$ is a
$\W_3$-manifold if and only if its associated pseudo-Riemannian
almost product manifold $(M,P,\tilde{g})$ is also a
$\W_3$-manifold.
\end{thm}
\begin{proof}
Let $(M,P,g)$ be a Riemannian $\W_3$-manifold. Using \eqref{22},
\eqref{23} and \eqref{sigma}, we obtain immediately
\begin{equation}\label{24}
\tilde{\nabla}_x y=\nabla_x y-\left(\nabla_x
P\right)Py-\left(\nabla_y P\right)Px,
\end{equation}
hence we have
\begin{equation}\label{25}
\left(\tilde\nabla_x P\right)y=-\left(\nabla_{Py}
P\right)Px-\left(\nabla_{x} P\right)y-\left(\nabla_{y} P\right)x.
\end{equation}
Since for a Riemannian $\W_3$-manifold $(M,P,g)$ the tensor
$\bar{N}$, defined by \eqref{6}, vanishes, then \eqref{25} has the
form
\begin{equation}\label{26}
\left(\tilde\nabla_x P\right)y=\left(\nabla_{Px} P\right)Py.
\end{equation}
According to properties $\tilde g(x,y)=g(x,Py)=\tilde g(Px,Py)$,
\eqref{2} and \eqref{26}, we get
\begin{equation}\label{27}
\tilde F(x,y,z)=-F(Px,y,z),
\end{equation}
from where
\begin{equation}\label{28}
\mathop{\s}_{x,y,z} \tilde F(x,y,z)=-\mathop{\s}_{x,y,z}
F(Px,y,z)=0.
\end{equation}
The equations \eqref{28} and \eqref{sigmaP} imply that
$(M,P,\tilde{g})$ is also a $\W_3$-manifold. Inversely, if
$(M,P,\tilde{g})$ is a $\W_3$-manifold, then, according to
\eqref{28} and \eqref{sigmaP}, we have that $(M,P,g)$ is also
$\W_3$-manifold.
\end{proof}

\begin{rem}
The condition \eqref{27} implies that a Riemannian $\W_3$-manifold
$(M,P,g)$ is a Riemannian $P$-manifold if and only if its
associated manifold $(M,P,\tilde{g})$ is a pseudo-Riemannian
$P$-manifold.
\end{rem}

\begin{thm}\label{thm:6}
Let $R$ be the curvature tensor of the Levi-Civita connection $\nabla$ for a
Riemannian $\W_3$-manifold $(M,P,g)$ and let $\tilde R$ be the curvature
tensor of the Levi-Civita connection $\tilde\nabla$ for its associated
manifold $(M,P,\tilde g)$. Then the following condition is valid
\begin{equation}\label{29}
\begin{split}
    \tilde R(x,y,z,w)=&R(x,y,z,Pw)-\left(\nabla_x F\right)(w,y,z)+\left(\nabla_y F\right)(w,x,z)\\[4pt]
        & +g\bigl(\left(\nabla_y P\right)z+\left(\nabla_z P\right)y,
            \left(\nabla_x P\right)Pw+\left(\nabla_w P\right)Px\bigr)       \\[4pt]
      & -g\bigl(\left(\nabla_x P\right)z+\left(\nabla_z P\right)x,
            \left(\nabla_y P\right)Pw+\left(\nabla_w P\right)Py\bigr).
 \end{split}
\end{equation}
\end{thm}
\begin{proof}
It is known that at the transformation $\nabla \rightarrow
\tilde{\nabla}$ it satisfies the following
\begin{equation}\label{30}
    \tilde R(x,y)z=R(x,y)z+Q(x,y)z,
\end{equation}
where
\begin{equation}\label{31}
    Q(x,y)z=\left(\nabla_x T\right)(y,z)-\left(\nabla_y T\right)(x,z)+T\left(x,T(y,z)\right)
    -T\left(y,T(x,z)\right),
\end{equation}
\begin{equation}\label{32}
    T(x,y)=-\left(\nabla_x P\right)Py-\left(\nabla_y P\right)Px.
\end{equation}
The corresponding tensor $\tilde R$ of type (0,4) is
\[
    \tilde R(x,y,z,w)=\tilde g\bigl(\tilde R(x,y)z,w\bigr).
\]
Hence, according to \eqref{30}, we have
\begin{equation}\label{34}
    \tilde R(x,y,z,w)=R(x,y,z,Pw)+ g\bigl(Q(x,y)z, Pw\bigr).
\end{equation}
Using \eqref{3}, \eqref{sigma}, \eqref{31} and \eqref{32} we
establish after some transformations the following
\begin{equation}\label{35}
\begin{split}
    g\bigl(Q(x,y)z, Pw\bigr)=&-\left(\nabla_x F\right)(w,y,z)+\left(\nabla_y F\right)(w,x,z)\\[4pt]
            & +g\bigl(\left(\nabla_y P\right)z+\left(\nabla_z P\right)y,
            \left(\nabla_x P\right)Pw+\left(\nabla_w P\right)Px\bigr)       \\[4pt]
      & -g\bigl(\left(\nabla_x P\right)z+\left(\nabla_z P\right)x,
            \left(\nabla_y P\right)Pw+\left(\nabla_w P\right)Py\bigr).
 \end{split}
\end{equation}
Then \eqref{34} and \eqref{35} imply \eqref{29}.
\end{proof}

\section{Invariant tensors of the transformation $\nabla \rightarrow \tilde{\nabla}$ }

Let $\nabla$ and $\tilde{\nabla}$ be the Levi-Civita connections
of a Riemannian $\W_3$-manifold $(M,P,g)$ and its associated
manifold $(M,P,\tilde g)$, respectively. An important problem is
the finding of invariant tensors of the transformation $\nabla
\rightarrow \tilde{\nabla}$.

We shall prove the following
\begin{thm}\label{thm7}
    Let $\nabla$ and $\tilde{\nabla}$ be the Levi-Civita connections of a Riemannian
    $\W_3$-manifold $(M,P,g)$ and its associated manifold $(M,P,\tilde g)$, respectively.
    Then at the transformation $\nabla \rightarrow \tilde{\nabla}$ the tensors
\begin{equation}\label{36}
    S(x,y)=\nabla_x y+\frac{1}{2}T(x,y),\quad L(x,y)z=R(x,y)z+\frac{1}{2}Q(x,y)z
\end{equation}
are invariant, where $T(x,y)$ and $Q(x,y)z$ are determined by \eqref{32} and \eqref{31}, respectively.
\end{thm}
\begin{proof}
Having in mind $\left(\tilde\nabla_x P\right)y=\tilde\nabla_x Py-P\tilde\nabla_x y$
and \eqref{24}, \eqref{32}, we get
\begin{equation}\label{37}
    \tilde T(x,y)=-\left(\tilde\nabla_x P\right)Py-\left(\tilde\nabla_y P\right)Px=-T(x,y).
\end{equation}
From \eqref{24} and \eqref{37} it follows
\[
    \tilde\nabla_x y+\frac{1}{2}\tilde T(x,y)=\nabla_x y+\frac{1}{2} T(x,y),
\]
from where we have
\[
    \tilde S(x,y)=\tilde\nabla_x y+\frac{1}{2}\tilde T(x,y)=S(x,y).
\]
This means that the tensor $S$ is invariant with respect to the transformation
$\nabla \rightarrow \tilde{\nabla}$.

At the transformation $\nabla \rightarrow \tilde{\nabla}$ the
tensor $Q$, determined by \eqref{31}, is transformed into the
tensor $\tilde Q$, determined by
\begin{equation}\label{40}
    \tilde Q(x,y)z=\left(\tilde \nabla_x \tilde T\right)(y,z)
    -\left(\tilde \nabla_y \tilde T\right)(x,z)+\tilde T\left(x,\tilde T(y,z)\right)
    -\tilde T\left(y,\tilde T(x,z)\right).
\end{equation}
Using \eqref{24}, \eqref{37} and \eqref{40}, we obtain
\begin{equation}\label{41}
    \tilde Q(x,y)z=-Q(x,y)z.
\end{equation}
Then \eqref{30} and \eqref{41} imply immediately
\[
    \tilde L(x,y)z=\tilde R(x,y)z+\frac{1}{2}\tilde Q(x,y)z=L(x,y)z,
\]
which means that the tensor $L$ is also an invariant one with
respect to the transformation $\nabla \rightarrow \tilde{\nabla}$.
\end{proof}

In the next theorems some characteristics of Riemannian $\W_3$-manifolds with vanishing
invariant tensors $S$ and $L$ are given.

\begin{thm}\label{thm:8}
    A Riemannian $\W_3$-manifold with zero tensor $S$ is a Riemannian $P$-manifold.
\end{thm}
\begin{proof}
    Let $(M,P,g)$ be a Riemannian $\W_3$-manifold with $S=0$. Then \eqref{36} and \eqref{32}  imply
\[
    \nabla_x y=\frac{1}{2}\bigl(\left(\nabla_x
P\right)Py+\left(\nabla_y P\right)Px\bigr),
\]
hence we obtain
\begin{equation}\label{44}
    \left(\nabla_y
    P\right)Px+\left(\nabla_{Py} P\right)x=0.
\end{equation}
In the last equation we substitute $x\leftrightarrow y$ and we have
\begin{equation}\label{45}
    \left(\nabla_x P\right)Py+\left(\nabla_{Px} P\right)y=0.
\end{equation}
Since \eqref{44} and \eqref{45}, the Nijenhuis tensor of $P$ has
the form
\[
N(x,y)=\left(\nabla_x P\right)Py+\left(\nabla_{Px} P\right)y
-\left(\nabla_y P\right)Px-\left(\nabla_{Py} P\right)x=0.
\]
It is known from \cite{StGr:connect} that $N(x,y)=0$ is a
characteristic condition $(M,P,g)$ to belong to the class $\W_1\oplus\W_2$.
Therefore  $(M,P,g)\in\left(\W_1\oplus\W_2\right)\cap\W_3=\W_0$, i.e.
$(M,P,g)$ is a Riemannian $P$-manifold.
\end{proof}

\begin{thm}\label{thm:9}
    Let $(M,P,g)$ be a Riemannian $\W_3$-manifold with zero tensor $L$. Then the following identity is valid
\begin{equation}\label{46}
\begin{split}
    2\bigl(&R(x,y,z,w)+R(x,Py,Pz,w)+\\[4pt]
    &+R(Px,Py,z,w)+R(Px,y,Pz,w)\bigr)=\\[4pt]
    &=2g\bigl(\left(\nabla_y P\right)Px+\left(\nabla_{Py} P\right)x,\left(\nabla_{Pw} P\right)z\bigr)\\[4pt]
    &\phantom{2}+g\bigl(\left(\nabla_y P\right)Pz+\left(\nabla_{Py} P\right)z,\left(\nabla_{Pw} P\right)x\bigr)\\[4pt]
    &\phantom{2}+g\bigl(\left(\nabla_z P\right)Px+\left(\nabla_{Pz} P\right)x,\left(\nabla_{Pw} P\right)y\bigr).
\end{split}
\end{equation}
\end{thm}
\begin{proof}
    Let $(M,P,g)$ be a Riemannian $\W_3$-manifold with $L=0$. According to \eqref{36} we have
\begin{equation}\label{47}
    R(x,y)z+\frac{1}{2}Q(x,y)z=0.
\end{equation}
We denote
\begin{equation}\label{48}
\begin{split}
B(x,y,z,w)=&-g\bigl(\left(\nabla_z P\right)y+\left(\nabla_y P\right)z,\left(\nabla_x P\right)w+\left(\nabla_{Pw} P\right)Px\bigr)\\[4pt]
          &+g\bigl(\left(\nabla_x P\right)z+\left(\nabla_z P\right)x,\left(\nabla_y P\right)w+\left(\nabla_{Pw} P\right)Py\bigr).\\[4pt]
\end{split}
\end{equation}
Then according to \eqref{47} and \eqref{35} we obtain
\begin{equation}\label{49}
    2R(x,y,z,w)=\left(\nabla_x F\right)(Pw,y,z)-\left(\nabla_y F\right)(Pw,x,z)+B(x,y,z,w).
\end{equation}
In \eqref{49} we substitute $Py$ and $Pz$ for $y$ and $z$, respectively, and then we add the obtained equation to \eqref{49}.
We get
\[
\begin{split}
    2\bigl(&R(x,Py,Pz,w)+R(x,y,z,w)\bigr)=\\[4pt]
    &=\left(\nabla_x F\right)(Pw,y,z)+\left(\nabla_x F\right)(Pw,Py,z)\\[4pt]
    &-\left(\nabla_y F\right)(Pw,x,z)-\left(\nabla_{Py} F\right)(Pw,x,Pz)\\[4pt]
    &+B(x,y,z,w)+B(x,Py,Pz,w).\\[4pt]
\end{split}
\]
Hence, according to \eqref{8}, we have
\begin{equation}\label{51}
\begin{split}
    2\bigl(&R(x,Py,Pz,w)+R(x,y,z,w)\bigr)=\\[4pt]
    &=A(x,Pw,Py,z)-A(Py,Pw,x,z)\\[4pt]
    &+B(x,y,z,w)+B(x,Py,Pz,w)\\[4pt]
    &-\left(\nabla_y F\right)(Pw,x,z)+\left(\nabla_{Py} F\right)(Pw,Px,z).\\[4pt]
\end{split}
\end{equation}
Now, in \eqref{51} we substitute $Px$ and $Py$ for $x$ and $y$, respectively, and
then we add the obtained equation to \eqref{51}.
Having in mind $\bar N(x,y)=0$, \eqref{6}, \eqref{8} and \eqref{48}, we get \eqref{46}.
\end{proof}

Using the last theorem we prove the following
\begin{cor}\label{thm:10}
        If $(M,P,g)$ is a Riemannian $\W_3$-manifold with zero tensor $L$, then
        \begin{equation}\label{thm10-1}
        \nP=-8\tau,
        \end{equation}
        \begin{equation}\label{thm10-2}
        \tau^{**}=5\tau.
        \end{equation}
\end{cor}
\begin{proof}
    If $(M,P,g)$ is a Riemannian $\W_3$-manifold with $L=0$, then we have \eqref{46}.
    Hence we obtain immediately
\begin{equation}\label{53}
\begin{split}
    2\bigl(&\rho(y,z)+\rho^*(Py,z)+\rho(Py,Pz)+\rho^*(y,Pz)\bigr)=\\[4pt]
    &=2g^{ij}g\bigl(\left(\nabla_y P\right)e_i+\left(\nabla_{Py} P\right)Pe_i,\left(\nabla_{Pe_j} P\right)Pz\bigr)\\[4pt]
    &\phantom{+}+g^{ij}g\bigl(\left(\nabla_z P\right)e_i+\left(\nabla_{Pz} P\right)Pe_i,\left(\nabla_{Pe_j} P\right)Py\bigr)\\[4pt]
    &\phantom{+}+g^{ij}g\bigl(\left(\nabla_y P\right)z+\left(\nabla_{Py} P\right)Pz,\left(\nabla_{Pe_i} P\right)Pe_j\bigr).\\[4pt]
\end{split}
\end{equation}
Since $(M,P,g)$ is a Riemannian $\W_3$-manifold, then $\theta=0$
and \eqref{nP} is valid. From \eqref{53}, using \eqref{Pg}, we
obtain $4(\tau+\tau^{**})=-3\nP$. Hence, according to \eqref{nP2},
we get \eqref{thm10-1} and \eqref{thm10-2}.
\end{proof}


\bigskip

\textit{Dimitar Mekerov\\
University of Plovdiv\\
Faculty of Mathematics and Informatics
\\
Department of Geometry\\
236 Bulgaria blvd.\\
Plovdiv 4003\\
Bulgaria
\\
e-mail: mircho@uni-plovdiv.bg}

\end{document}